\documentclass[12pt]{amsart}%
\usepackage{color}
\usepackage{amsmath}
\usepackage{amssymb}
\usepackage{amsfonts}
\usepackage[latin1]{inputenc}
\usepackage{graphicx}%
\usepackage{cite}
\providecommand{\U}[1]{\protect\rule{.1in}{.1in}}

\DeclareMathSymbol{\subsetneqq}{\mathbin}{AMSb}{36}

\theoremstyle{plain}
\numberwithin{equation}{section}
\newtheorem{theorem}{Theorem}[section]

\newtheorem{lemma}{Lemma}[section]

\newtheorem{remark}{Remark}[section]
\newtheorem{notation}{Notation}
\input{tcilatex}
\begin{document}
\title[A simple proof of the Karush-Kuhn-Tucker theorem]
{A simple proof of the Karush-Kuhn-Tucker theorem with finite number of equality and inequality constraints}%
\author{Ramzi May}%
\address{Department of Mathematics and Statistics, College of Science, King Faisal University, Al-Ahsa, Kingdom of Saudi Arabia}
\email{rmay@kfu.edu.sa}
\vskip 0.2cm
\date{ June 21, 2020}

\begin{abstract} We provide a simple and short proof of the
Karush-Kuhn-Tucker theorem with finite number of equality and inequality
constraints. The proof relies on an elementary linear algebra lemma and the
local inverse theorem.
\end{abstract}
\maketitle

\section{Introduction}

Let $X$ be a normed real linear space. We denote by $X^{^{\prime }}$ the
space of linear mapping from $X$ to $\mathbb{R}$ and by $X^{\ast }$ the
dual space of $X$ i.e. the space of linear and continuous mapping from $X$
to $\mathbb{R}$. The infinite dimensional version of the famous
Karush-Kuhn-Tucker theorem with finite number of equality and inequality
constraints reads as follows.

\begin{theorem}
\label{Th}
let $\Omega $ be an open space of $X$ and $\{f_{i}: 0\leq i\leq n+m \}$ a family of continuously differentiable functions from $\Omega $
to $\mathbb{R}$ where $n\in \mathbb{N},m\in \mathbb{N}\cup \{0\}$.  Let $x^{\ast }$
be a solution of the constraint minimization problem%
\[
\left\{
\begin{array}{ll}
\text{minimize} & f_{0}(x) \\
\text{subject to}
& f_{i}(x)=0,1\leq i\leq n, \\
& f_{i}(x)\leq 0,n+1\leq i\leq n+m,%
\end{array}%
\right.
\]%
such that the family $\{f_{i}^{\prime }(x^{\ast }): i\in J(x^{\ast })\}$ is linearly independent in $X^{\ast },$ where
\[
J(x^{\ast })=\{i: 1\leq i\leq n+m~\text{and}~ f_{i}(x^{\ast })=0\}.
\]%
Then there exist $(\lambda _{1},\cdots ,\lambda _{n})\in \mathbb{R}^{n}$ and
$(\mu _{1},\cdots ,\mu _{m})\in ([0,+\infty \lbrack )^{m}$ such that
\[
f_{0}^{\prime }(x^{\ast })+\sum_{i=1}^{n}\lambda _{i}f_{i}^{\prime }(x^{\ast
})+\sum_{j=1}^{m}\mu _{j}f_{j+n}^{\prime }(x^{\ast })=0
\]%
and%
\[
\mu _{j}f_{j+n}(x^{\ast })=0,~\forall 1\leq j\leq m.
\]
\end{theorem}

This famous theorem is a natural extension of the classical
Lagrange multipliers theorem to the case of the minimization problem with
finite number of equality and inequality constraints. Its finite dimensional
version has been originally derived independently by Karush \cite{K} and Kuhn and
Tucker \cite{KT}. Since there, different proofs of the generalization of the Karush,
Kuhn and Tucker theorem (KKT) to the infinite dimensional setting have been
provided in many works (see, for instance, \cite{BS,CZ,GT,R} and references therein). In three recent papers \cite{BTW0,BT1,BTW2}, Brezhneva, Tretyakov, and Wright have given some elementary and different proofs of the KKT Theorem respectively with equality constraints, inequality constraints and linear equality, and nonlinear inequality constraints. In this short note, inspired essentially by the paper \cite{BTW2}, we give a new,  detailed
and simple proof of the KKT theorem with finite number of mixed equality and
inequality constraints. Our proof relies essentially on a very simple but
powerful lemma from linear algebra and the classical local inverse theorem
in the finite dimensional setting.

\section{Proof of the Karush, Kuhn and Tucker Theorem}
Before starting the proof of Theorem \ref{Th}, we introduce the following simple notations.
\begin{notation}
Let $N\in \mathbb{N} $.
\begin{enumerate}
\item The canonical basis of $\mathbb{R}^N$ is the vector family $\{e_1,\cdots,e_N\}$ defined be: $e_1=(1,0,\cdots,0), e_2=(0,1,0,\cdots,0), \cdots, e_N=(0,\cdots,0,1)$.
\item $B_{\mathbb{R}^N}(0,r)$ is the open ball of $\mathbb{R}^N$ with center $0$ and radius $r>0$.
\item $I_N$ is the unity matrix of size $(N,N)$.
\item  $Id_{\mathbb{R}^{N}}$ is the identity mapping from $\mathbb{R}^{N}$ into itself.
\end{enumerate}
\end{notation}

Next, we prove the following elementary linear algebra lemma.

\begin{lemma}
\label{le}
Let $\{T_{i}: 1\leq i\leq n\}$ be a finite family of a linear
independent elements of $X^{\prime }.$ Then there exists a family $%
\{v_{i}: 1\leq i\leq n\}$ of elements of $X$ such that
\begin{equation}
T_{i}(v_{j})=\delta _{ij},~\forall 1\leq i,j\leq n  \label{lin}
\end{equation}%
where $\delta _{ij}$ is the Kronecker's symbol.
\end{lemma}

\begin{remark}
The family $\{v_{i}: 1\leq i\leq n\}$ will be called a quasi
primal basis of $X$ associated to the family $\{T_{i}: 1\leq i\leq n\}.$
\end{remark}

\begin{proof}
Define the linear mapping $T:X\rightarrow \mathbb{R}^{n},~T(v)=(T_{1}(v),%
\cdots ,T_{n}(v)).$ Let us prove that $T$ is
not surjective. Suppose that this is not true; then there exists a vector $\alpha =(\alpha _{1},\cdots
,\alpha _{n})\in \mathbb{R}^{n}\backslash \{0\}$ orthogonal in $\mathbb{R}%
^{n}$ (with respect to the usual inner product) to the linear subspace $%
T(X) $; which implies that for every $v\in X,$%
\[
\alpha _{1}T_{1}(v)+\cdots +\alpha _{n}T_{n}(v)=0.
\]%
This contradicts the assumption on the family $\{T_{i}:~1\leq i\leq n\}.$
Therefore we conclude that the mapping $T$ is surjective. By consequence, for every $%
1\leq j\leq n,$ there exists $v_{j}\in X$ such that $%
T(v_{j})=e_{j}$, where $\{e_{1},\cdots ,e_{n}\}$ is the canonical basis of $%
\mathbb{R}^{n}.$ Clearly, the family $\{v_{i}: 1\leq i\leq n\}$
satisfies (\ref{lin}).
\end{proof}

Now we are ready to prove the KKT Theorem.
\begin{proof}
Let us first notice that up to replace $\Omega $ by the open subset
\[
\Omega ^{\ast }=\{x\in \Omega :f_{i}(x)<0,~ \forall n\leq i\leq n+m ~ \text{and} ~ i \notin J(x^{\ast })\}
\]%
and to set $\mu _{j-n}=0$ for every $j\notin J(x^{\ast
}),$ we can assume without loss of generality that $J(x^{\ast })=\{i: 1\leq i\leq n+m\}
.$ Now we will first prove that the family $\{f_{i}^{\prime }(x^{\ast
}): 0\leq i\leq n+m\}$ is linearly dependant in $X^{\ast }.$ We
argue by contradiction. According to Lemma \ref{le}, there exits a quasi primal
basis $\{v_{i}: 0\leq i\leq n+m\}$ of $X$ associated to the
family $\{f_{i}^{\prime }(x^{\ast }): 0\leq i\leq n+m\}.$ Since $%
x^{\ast }$ belongs to the open subset $\Omega $, there exists a real number $%
r_{0}>0$ such that the mapping defined for every $t=(t_{i})_{0\leq i\leq n+m}$ in $B_{\mathbb{R}^{m+n+1}}(0,r_{0})$  by
\[
\Phi (t)=(f_{0}(\sigma (t)),\cdots ,f_{m+n}(\sigma (t))),
\]%
where
\[
\sigma (t)=x^{\ast }+\sum_{i=0}^{m+n}t_{i}v_{i},
\]%
is continuously differentiable and its Jacobian matrix at $t=0$ is
\[
J_{\Phi }(0)=\left[ f_{i}^{\prime }(x^{\ast })(v_{j})\right] _{0\leq i,j\leq m+n}=I_{m+n+1}.
\]%
Therefore, $\Phi ^{\prime }(0)=Id_{\mathbb{R}^{m+n+1}}$; hence by applying
the local inverse theorem, we deduce the existence of a real number $%
r_{1}\in ]0,r_{0}]$ such that $\Phi $ is a $C^{1}$ diffeomorphism from $%
U_{1}\equiv B_{\mathbb{R}^{m+n+1}}(0,r_{1})$ to an open neighbourhood $V_{1}$
of $\Phi (0)=(f_{0}(x^{\ast }),0,\cdots ,0)$ in $\mathbb{R}^{m+n+1}.$ For $%
\nu >0$ small enough, the vector $y_{\nu }\equiv (f_{0}(x^{\ast })-\nu
,0,\cdots ,0)$ belongs to $V_{1};$ let $t_{\nu }=\Phi ^{-1}(y_{\nu }).$ It
is clear that the vector $x_{\nu }=\sigma (t_{\nu })$ belongs to $\Omega $
and satisfies
\begin{eqnarray*}
f_{0}(x_{\nu }) &=&f_{0}(x^{\ast })-\nu , \\
f_{i}(x_{\nu }) &=&0,~\forall 1\leq i\leq n+m,
\end{eqnarray*}%
which contradicts the definition of $x^{\ast }.$ Thus, the
family $\{f_{i}^{\prime }(x^{\ast }): 0\leq i\leq n+m\}$ is
linearly dependant in $X^{\ast }.$ On the other hand, since $%
\{f_{i}^{\prime }(x^{\ast }): 1\leq i\leq n+m\}$ is linearly
independent in $X^{\ast },$ we infer the existence of $(\lambda _{1,\cdots
,}\lambda _{n},\mu _{1,}\cdots ,\mu _{m})\in \mathbb{R}^{m+n}$ such that
\begin{equation}
f_{0}^{\prime }(x^{\ast })+\sum_{i=1}^{n}\lambda _{i}f_{i}^{\prime }(x^{\ast
})+\sum_{j=1}^{m}\mu _{j}f_{j+n}^{\prime }(x^{\ast })=0.  \label{A1}
\end{equation}%
It remains to prove that $\mu _{j}\geq 0$ for every $1\leq j \leq m$
According to Lemma \ref{le}, there exists $\{w_{1},\cdots ,w_{m+n}\}$ a quasi primal
basis of $X$ associated to the family $\{f_{i}^{\prime }(x^{\ast }): 1\leq i\leq n+m\}.$ Proceeding as previously, we
deduce that there exists $r>0$ and a neighbourhood $V$ of $0$ in $\mathbb{R}%
^{m+n}$ such that the mapping $\varphi :B_{\mathbb{R}^{m+n}}(0,r)\rightarrow
V$ defined
\[
\varphi (t)=(f_{1}(s(t)),\cdots ,f_{m+n}(s(t)),
\]%
where%
\[
s(t=(t_{i})_{1\leq i\leq m+n})=x^{\ast }+\sum_{i=1}^{m+n}t_{i}w_{i},
\]%
is a $C^{1}$ diffeomorphism. Let $1\leq j_{0}\leq m$ be a fixed
integer. Since $V$ is an open neighbourhood of $0$ in $\mathbb{R}^{m+n},$
there exists $\varepsilon _{0}>0$ such that for every $\varepsilon \in
]-\varepsilon _{0},\varepsilon _{0}[,~-\varepsilon e_{j_{0}+n}\in V,$
where $(e_{1},\cdots ,e_{m+n})$ is the canonical basis of $%
\mathbb{R}^{m+n}.$ Hence, for every $\varepsilon \in ]-\varepsilon
_{0},\varepsilon _{0}[,$ the vector
\[
\tilde{x}(\varepsilon )=s(\varphi ^{-1}(-\varepsilon e_{j_{0}+n}))
\]%
belongs to $\Omega $ and satisfies $f_{j_{0}+n}(\tilde{x}(\varepsilon
))=-\varepsilon $ and $f_{i}(\tilde{x}(\varepsilon ))=0$ for every $i\in \{1,\cdots,m+n\}\backslash j_{0}+n.$ Hence, for every $\varepsilon
\in ]0,\varepsilon _{0}[,$%
\[
\frac{f_{0}(\tilde{x}(\varepsilon ))-f_{0}(\tilde{x}(0))}{\varepsilon }=%
\frac{f_{0}(\tilde{x}(\varepsilon ))-f_{0}(x^{\ast })}{\varepsilon }\geq 0.
\]%
Letting $\varepsilon \rightarrow 0,$ we obtain%
\begin{equation}
f_{0}^{\prime }(x^{\ast })(\frac{d\tilde{x}}{d\varepsilon }(0))\geq 0.
\label{A2}
\end{equation}%
For every $\varepsilon \in ]-\varepsilon _{0},\varepsilon _{0}[,$ define
\[
\tilde{t}(\varepsilon )=(\tilde{t}_{1}(\varepsilon ),\cdots ,\tilde{t}%
_{m+n}(\varepsilon ))=\varphi ^{-1}(-\varepsilon e_{j_{0}+n}).
\]%
First, since $\varphi (\tilde{t}(\varepsilon ))=-\varepsilon e%
_{j_{0}+n}$ and $\varphi ^{\prime }(0)=Id_{\mathbb{R}^{m+n}},$ we have $%
\frac{d\tilde{t}}{d\varepsilon }(0)=-e_{j_{0}+n}.$ Using now the
fact that%
\[
\tilde{x}(\varepsilon )=s(\tilde{t}(\varepsilon ))=x^{\ast }+\sum_{i=1}^{m+n}%
\tilde{t}_{i}(\varepsilon )w_{i},
\]%
we deduce that%
\begin{equation}
\frac{d\tilde{x}}{d\varepsilon }(0)=-w_{j_{0}+n}.  \label{A3}
\end{equation}%
Finally, combining (\ref{A1}),(\ref{A2}), and (\ref{A3}) yields%
\[
\mu _{j_{0}}=-f_{0}^{\prime }(x^{\ast })(w_{j_{0}+n})\geq 0,
\]%
which completes the proof of the theorem.
\end{proof}

\end{document}